\documentclass[11pt]{article}
\usepackage{hfstyle}
\usepackage[utf8]{inputenc}

\usepackage{amsmath,graphicx, color}

\DeclareFixedFont{\ttb}{T1}{txtt}{bx}{n}{12} % for bold
\DeclareFixedFont{\ttm}{T1}{txtt}{m}{n}{12}  % for normal

\usepackage{color}
\definecolor{deepblue}{rgb}{0,0,0.5}
\definecolor{deepred}{rgb}{0.6,0,0}
\definecolor{deepgreen}{rgb}{0,0.5,0}

\usepackage{listings}

\newcommand\pythonstyle{\lstset{
language=Python,
basicstyle=\ttm,
morekeywords={self},              % Add keywords here
keywordstyle=\ttb\color{deepblue},
emph={MyClass,__init__},          % Custom highlighting
emphstyle=\ttb\color{deepred},    % Custom highlighting style
stringstyle=\color{deepgreen},
frame=tb,                         % Any extra options here
showstringspaces=false
}}

\lstnewenvironment{python}[1][]
{
\pythonstyle
\lstset{#1}
}
{}
\begin{document}
\title{Mathematics in the liturgical books of the Catholic Church: phases of the ecclesiastical moon}
\author{Henryk Fuk\'s  
      \oneaddress{
         Department of Mathematics and Statistics\\
         Brock University,  St. Catharines, Canada\\
         \email{hfuks@brocku.ca}
       }
   }

\Abstract{
We use contemporary mathematical  notation
to describe the  method for determining the age of the ecclesiastical moon as mandated by pope Gregory XIII and elaborated in the book of Christopher Clavius \emph{Romani calendarii explicatio}. The algorithm is first introduced by 
using the tabular method employed  by liturgical books such as  the Roman Missal, Breviary  and Martyrology. Then we construct the recurrence equation for the epacts, derive its solution, and
give a simple expression for the age of the moon 
on a given day of the year.  We also consider the problems which can occur at the transition from December 31 to January 1 of the next year, when
there could be a ``jump'' in moon's age (\emph{saltus lunae}) in years when 
epact corrections are applied. We propose a simple solution
which fixes these problems. A summary of the formulae and listing of the 
implementation of relevant functions in Python is  provided in the last section.
}

\maketitle

\section{Gregorian calendar and the ecclesiastical moon}
The Julian calendar, introduced by Julius Caesar in 46 BC 
\cite{Feeney2007},
has been used in the Western world for the next sixteen
 centuries. 
The average length of the year in the Julian
calendar was 365.25 days, and this was
slightly more than  the actual duration of the mean solar year (currently  365.24217 mean solar days). As a result,
the Julian calendar was gaining a day every 128 years, 
thus falling behind the astronomical year more and more.
In the 16-the century, the astronomical vernal equinox was occurring 10 days before its nominal date of March 21.

The problem was known since Middle Ages \cite{nothaft2018scandalous}, but
toward the end of the 15th century the progress in astronomy 
intensified the  calls for the calendar reform.
The reform was eventually undertaken during the reign of Pope Gregory XIII (1572–1585). 
It is not the purpose of this paper to  discuss the details of the history of the reform and all preceding
events, but the  interested reader can consult, for example, 
reference \cite{coyne1983gregorian}, where  numerous other works on the subject are sourced.

On February 24,  1582,  Pope Gregory XIII issued the bull
\emph{Inter gravissimas} introducing the new calendar, later called Gregorian. The new calendar reform corrected the length of the year by making century years not divisible by 400
non-leap years. This resulted in the average year
of 365.2425 days, which was much closer to the actual solar year.
Furthermore, the reform deleted 10 days from the calendar so that in the year of the reform,  October 4 was followed by October 15, bringing the date of the vernal equinox back to 
March 21.

While the above aspects of the reform are widely known, 
it is perhaps not so widely known that 
the Gregorian reform also corrected the way moon phases are
calculated. Knowledge of the moon phase on a given day 
is needed for liturgical purposes. First of all, 
it is needed to determine the
date of Easter,  which falls on the first Sunday after the 
the first full moon on or after  March 21.
Additionally, in many Catholic religious communities, the current phase of the moon
is pronounced at the beginning of daily reading
from the Roman Martyrology book.

The phase of the moon on a given day is usually expressed by
the \emph{moon age}, defined as the time elapsed since the last new moon. It is important to emphasize that the moon age used by the Church for liturgical purposes is not the age of the real astronomical moon. It is rather the age of the imaginary
moon, known as \emph{ecclesiastical moon}, which is determined
by application of a prescribed algorithm.
The crucial part of this algorithm, which we will discuss
below, is computation of the epacts. Epact of the year (lat. \emph{epacta}) is the the age of the ecclesiastical moon on January 1, minus~1.

The new  algorithm for computing epacts and thus determining the 
age of the ecclesiastical moon was
originally proposed by Aloysius Lilius, expanded upon by Christopher Clavius, and then mandated by Pope Gregory XIII
to be used in the Catholic Church.
Although the phase of the ecclesiastical moon computed with this algorithm  closely follows the actual phase of the physical  moon, it deviates from it slightly, partially due to the fact that the ecclesiastical moon's age is always expressed in whole days while the real moon's phase is continuously changing, and partially because 
Clavius designed it that way, for reasons
which have nothing to do with astronomy\footnote{Clavius wanted to make sure that the Easter Sunday does not coincide with the Jewish Pascha and that the new moon does not occur on the same date in two different years of a single Metonic cycle of 19 years.}.

The principal goal of this article is to present a concise formulae  for determining the age of the ecclesiastical moon using contemporary mathematical notation. While
parts of such formulae can be found in various works
in one form or another, I am not aware of any
single reference giving such formulae in full and discussing them in a comprehensive way. The second goal is to propose
a corrected version of the ecclesiastical moon age, 
fixing some problems which occasionally occur between December 31 and January 1 of the next year.

Before we get to the mathematical part, we will first
describe the method of determining the age of the 
ecclesiastical moon using liturgical books, or more
precisely, using tables included in various
liturgical books. In subsequent sections we will give the mathematical background and the justification of this method.

\section{Determining the age of the moon using liturgical books}
Details of the algorithm for determining the age of the moon
on a given day are described in the introductory sections of
the Roman Missal (\emph{Missale Romanum}) or
the Breviary (\emph{Breviarum Romanum}). It does not really matter which edition one uses, as all editions published after 
1582 describe it in the same way.

\begin{figure}[htb]
\begin{center}
 \includegraphics[width=12cm]{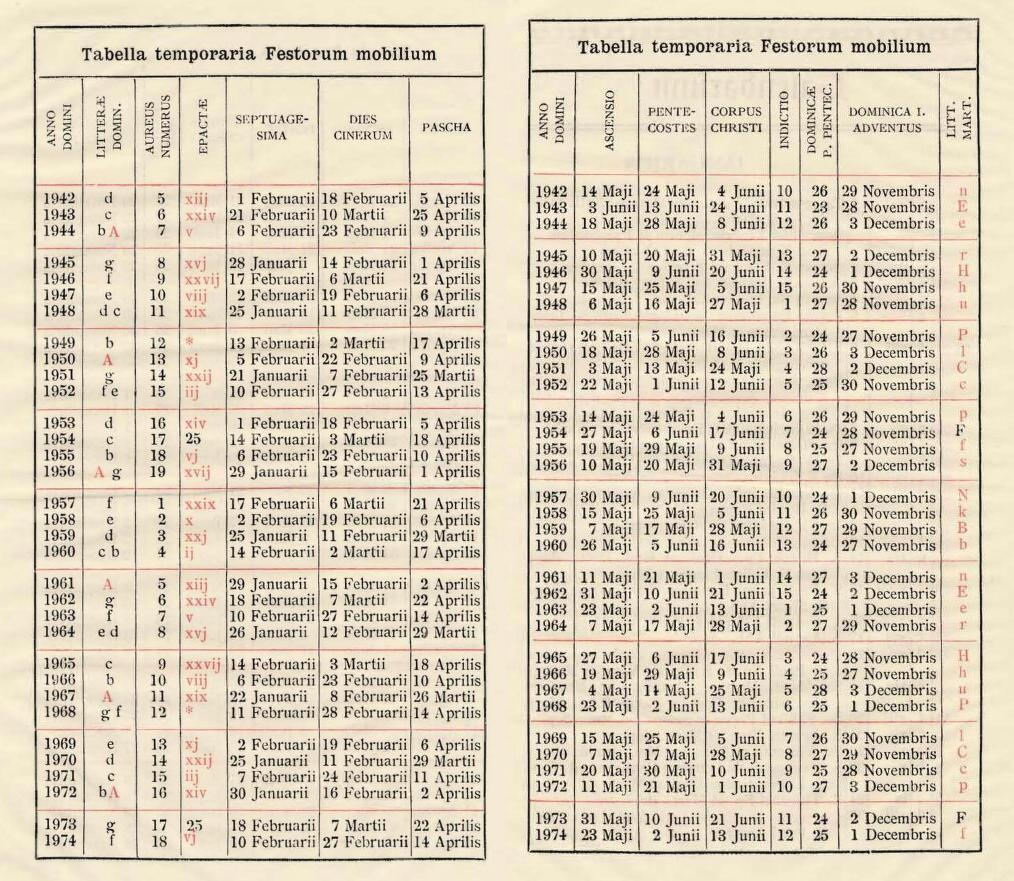}
\end{center}
\caption{\emph{Tabella temporaria} or temporal table reproduced from  older \emph{Breviarium Romanum}~\cite{brev1942}.}\label{tabtemp}
\end{figure}

In order to determine the age of the moon one needs to 
figure out three items, namely  the golden number, the epact, and
 the martyrology letter of the year. In \emph{Missale Romanum} 
or \emph{Breviarium Romanum}
these are usually given in a form of the table called \emph{Tabella temporaria}, constructed for years covering several decades following the publication date, 
like the one shown in Figure~\ref{tabtemp}.

We start with the golden number. 
 The golden number is the remainder from division of the year by 19, plus 1. It is used to label consecutive years of the 
 Metonic cycle of 19 years, after which the lunar phases recur at the same time of the year. In \emph{Tabella temporaria}
the golden number  is taken from the 
 column labeled \emph{aureus numerus}. For example, for 1945,  we see \emph{aureus numerus} 8.

The epact is the next item we need. As mentioned earlier,
 the epact is the  age of the ecclesiastical moon on January 1, minus 1.  For the real physical moon, the average 
 interval between phases is 29.53059 days, but 
can vary up to several hours from that number.  
 For the ecclesiastical
 moon, as we already mentioned, the age of the moon is measured
 in whole days, with either 29 or 30 days interval between phases.
For this reason,  the epact can take only integer values from 0 to 29.
The epact is given in the column \emph{Epacta}, where, for example, for 1945 we read ``xvj'', which in Arabic numerals corresponds to 16. It is customary 
to denote the epacts from 1 to 29 by Roman numerals i, ii, ii, iv, \ldots xxix. In older books if the last digit is ``i'', it is usually printed as ``j'', hence 
xvj instead of xvi.
If the epact is 0 (the last day of the cycle), it is denoted by~*.
Occasionally the epact will be listed in the table as 25 (in Arabic numerals), as it is the case for 1973 in  Figure~\ref{tabtemp}. This is different from ``xxv'', and we will see later why. The special epact 25 is used in years when the golden number is greater or equal to 12. For years with a golden number less than 12 the regular Roman number xxv is used. 

\begin{figure}
\begin{center}
\includegraphics[width=12cm]{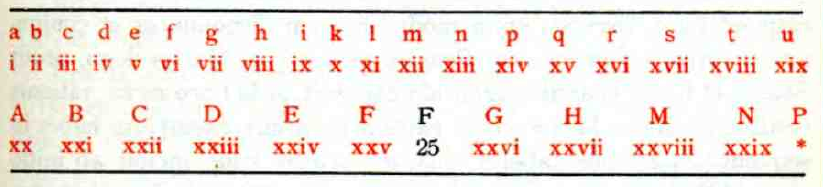}
\end{center}
\caption{Correspondence between martyrology  letters (top rows) and epact numbers (bottom rows), reproduced from
\emph{Martyrologium Romanum} \cite{mart1948}.
}\label{martlet}
\end{figure}

The third item we need from the \emph{Tabella temporaria} is the martyrology letter, \emph{littera martyrologii}, given in the last column. For 1945, this letter is ``r''.
Th martyrology letter is somewhat superficial, as it is a letter which simply corresponds to the epact. Epact i is represented by a, ii by b, iii by c, etc., as shown in Figure~\ref{martlet}. The special epact 25 is represented by F with a different colour than the remaining letters (in
Figure~\ref{martlet} it is black).
%%%%%%%%%%%%%%

Having the golden number, the epact,  and the martyrology letter taken from the temporal table, we can now come back to the determination of the age of the moon.
By adding 1 to the epact, we obtain
the age of the moon on January 1. Thus for 1945, on January 1 the moon age will be 8+1=9. If the epact is *, then the age of the moon on January 1 would be 1, i.e., the first day of the new moon. This would be the case for 1968, as we can see in Figure~\ref{tabtemp}.

\begin{figure}[htb]
\begin{center}
\includegraphics[width=12cm]{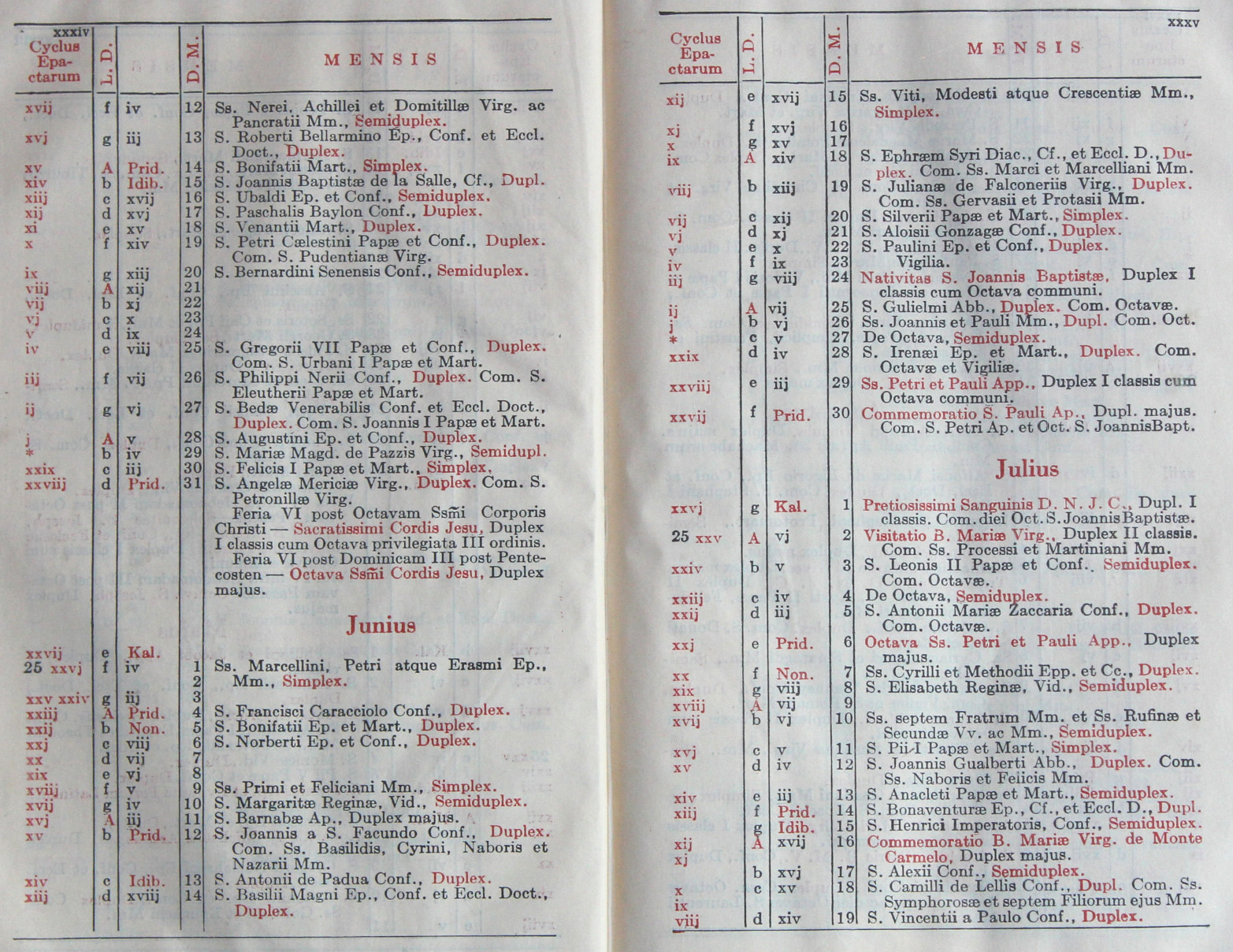} 
\end{center}
\caption{Fragment of \emph{Calendarium} from \emph{Missale Romanum} \cite{mr1942}.}\label{calendarium}
\end{figure}

What about the other days then, beyond January 1? There
are two ways of determining the age of the moon on other days. If all we have is the \emph{Breviarium Romanum} or the 
\emph{Missale Romanum}, we can consult \emph{Calendarium} table which is included at the beginning of these books, following the introductory section titled \emph{De Anno et Ejus Partibus}, ``About the year and its parts''. An  example fragment of the \emph{Calendarium}  is shown in Figure~\ref{calendarium}.
In the  \emph{Calendarium} all days of the year are listed
in consecutive order, and in the first column we see ``Cyclus Epactorum'', which can be used to find days of the new moon throughout  the year. If we know the epact of the year, then each day where that epact appears in the first column of the \emph{Calendarium} is the day of the new moon. For example, for 1945, we found that the epact of the year is xvj. In Figure~\ref{calendarium} we see fragments of \emph{Calendarium} corresponding to part of May, June and part of July. The epact xvj can be found on 
May 13, June 11, and July 11. All these days correspond to the first day of the moon. If we want to determine, for example, the age of the moon on July 15, 1945, we need to count forward to July 15 from the nearest new moon preceding July 15, namely July 11. Since July 11 is the first day of the moon, 
July 12 is the second day, and so forth. It is easy to see that July 15 will, therefore,  be the fifth day of the moon.

There is also another,  somewhat more convenient way of determining the age of the moon, using a different liturgical book, namely \emph{Martyrologium Romanum},  the Roman Martyrology. \emph{Martyrologium} lists martyrs and saints
commemorated on each day of the year, and entries for a given
day are to be read on the previous day  at the canonical Hour of Prime by people praying the Breviary.

\begin{figure}
\begin{center}
\includegraphics[width=12cm]{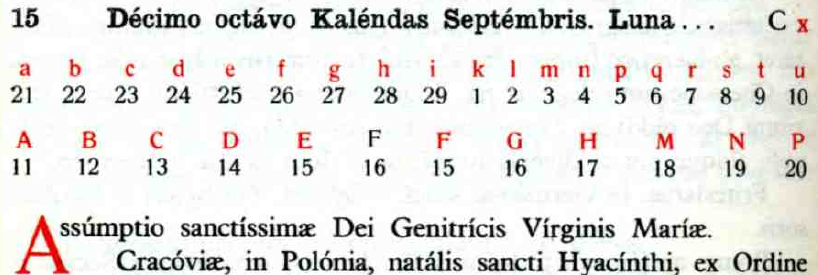} 
\end{center}
\caption{Fragment of the pre-concilliar \emph{Martyrologium Romanum}
\cite{mart1948} for August 15.
}\label{samplemartday}
\end{figure}

Rubrics of the Breviary stipulate that when the Martyrology is read publicly on a given day, the age of the moon can be optionally pronounced at the beginning. To help to determine the age of the moon quickly, special ``lunar table'' is provided for each day of the year. Example is shown in
Figure~\ref{samplemartday}. This figure shows an excerpt from
the page of \emph{Martyrologium} corresponding to  August 15, the feast of the Assumption. In order to determine the day 
of the moon one simply reads from the aforementioned table the number located below the letter corresponding to the current \emph{littera martyrologii}, martyrology letter. For 1945 the letter is ``r'', and under ``r'' we read 7, therefore August 15, 1945 was the seventh day of the moon.

\begin{figure}
\begin{center}
\includegraphics[width=12cm]{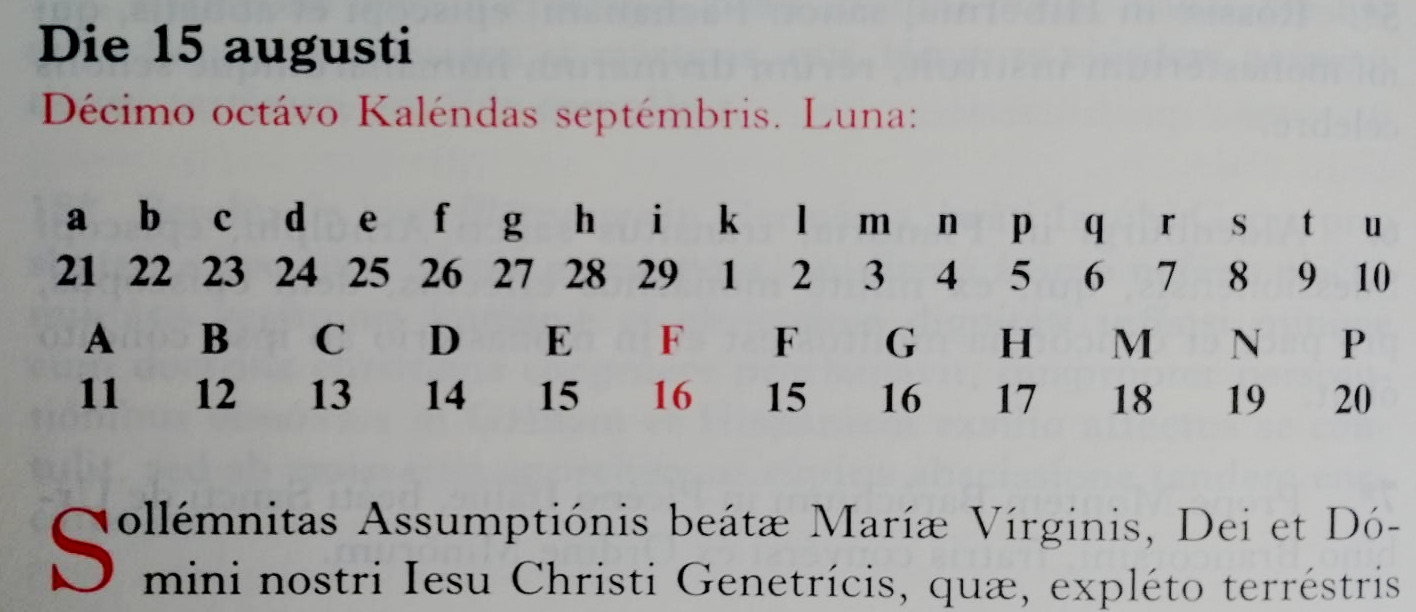} 
\end{center}
\caption{Fragment of the new (2004) \emph{Martyrologium Romanum}  \cite{mart2004} for August 15. Note the inverted color scheme.}\label{samplemartday-new}
\end{figure}

One detail is worth noticing at this point. In all pre-Vatican II editions of the Roman Martyrology, the martyrology letters were always printed in red ink, except the special case corresponding to Arabic epact 25, when  black F was used, as in Figure~\ref{samplemartday}. This was the case in all editions of the Martyrology starting from the first one published in 1583 till the last preconciliar
edition published in 1956. Unfortunately, for some incomprehensible reason,
in the post-conciliar edition of 2004, this colouring scheme has been reserved, so that all letters are black with the exception of F corresponding to the special epact  25,  which in the new edition is red, as shown in Figure~\ref{samplemartday-new}.
This can be a source of confusion and mistakes for Martyrology users in years with the epact 25. For example, in 1954, as
we can see in Figure~\ref{tabtemp}, the epact is 25,
and the martyrology letter is F. In the old \emph{Martyrologium}
(Figure \ref{samplemartday}) we will find the corresponding moon age under the black F, but in the new one (Figure \ref{samplemartday-new})  under the red F.
It is, therefore, best to remember that the letter F  colored differently than other letters corresponds to the special epact 25, while the letter F colored as other letters corresponds to xxv.
 %%%%%%%%%%%%%%%%%
%%%%%%%%%%%%%%%%%

\section{The recurrence equation for epacts}
As promised, we will now describe a more ``algebraic''
method for determining the age of the moon.
We will use $E_y$ to denote the epact of year $y$.
In the algorithm for the age of the moon it is assumed that $E_y$ can take 
values from 0 to 29 inclusively\footnote{In 
liturgical books and in the book of Clavius the value 0 is not used, the number 30 is used instead, represented by *. This does not really matter as we will always take the epact modulo 30.}.
If the epact of year $y$ is $E_y$, then the next year's epact
can be calculated recursively, following the algorithm given by Aloysius Lilius and implemented by the Jesuit mathematician Christopher Clavius \cite{Clavius1603}. 
The next year's epact is obtained 
by adding 11 days to the previous epact plus some corrections, and taking the result modulo 30.
We add 11 days because in the solar year there are 12 full moon cycles leaving (approximately) extra 11 days,
$365-12 \times 29.5=11$.
 For example, if the first new moon in a given year occurs on  January 12, then the next year's first new moon will be roughly 11 days earlier, on January 1. 
  We say ``roughly'' because some corrections are 
needed to trace the phase of the actual moon more precisely.

Here is the formal way to describe the recursion. 
We start by taking the epact of year 1582 to be 26, and then,
 as mentioned,
we compute the next year's epact by adding to the previous epact 11 days and  some corrections, modulo 30. The recurrence equation 
satisfied by $E_y$ is then
\begin{align}\label{epactrecursion}
&E_{1582}=26, \nonumber \\
&E_y=E_{y-1}+11 + M_y- S_y+L_y \mod 30.
\end{align}
In the above,  $M_y$, $S_y$ and $L_y$ are corrections stipulated in Clavius' book \cite{Clavius1603}. They are always zero except:
\begin{itemize}
\item $M_y=1$ for years $1596, 1615, 1634, 1653, 1672, 1691, 1710, 1729, 1748, \ldots$
\item $S_y=1$ for years 
$1700,1800,1900,2100,2200,2300,2500,2600,2700,\ldots$.
\item $L_y=1$ for years $1800, 2100, 2400, 2700, 3000, 3300, 3600, 3900,4300, \ldots$.
\end{itemize}
These corrections are described in detail as follows.
%Let $$C_y=\lfloor y/100 \rfloor +1$$ be the century number of %year $y$, and
Let
$$G_n=y \mod 19 +1$$ be the so-called  \emph{golden number} of $y$.
It is important in the epact's calculations because after 19 years (also known as Metonic cycle), the moon phases occurs on the same date, and this follows from the fact that there are almost exactly 235 synodic months in 19 tropical years.  
By synodic month we mean the is the average period of the Moon's orbit with respect to the line joining the Sun and Earth,
which is 29 d 12 h 44 min and 2.9 s. The tropical
year is the time that the Sun takes to return to the same position in the sky, equal to 365.24217 mean solar days.

$M_y$, which we will call \emph{Metonic cycle completion}, is equal to $1$ only when the golden number of the year is 1.
This correction is needed due to the fact that after the full Metonic cycle of 19 years, if we increase the epact by 11 days  each time, the total increase would be $19 \times 11 = 209$, which is not  divisible by 30 ($209= 6 \times 30 + 29)$. Adding 1 every 19 years is thus needed for the cycle of epacts to repeat itself. Years which are multiples of 19 (corresponding to $G_y=1)$ are the  correction years, therefore,
\begin{equation}\label{metoncor}
M_y = \begin{cases}
  1  & \text{if $G_y=1$,} \\
  0 &  \text{ otherwise.}
\end{cases}
\end{equation}
$S_y$ is called \emph{solar equation correction}, and it is
equal to 1 only for century years not divisible by 400, that is,
for $y\in \{1700, 1800, 1900, 2100, \ldots\}$. 
\begin{equation}\label{solcor}
S_y = \begin{cases}
  1  & \text{if $y$ mod $100=0$ and $y$ mod $400 \neq 0$,} \\
  0 &  \text{ otherwise.}
\end{cases}
\end{equation}
This accounts for the feature of the Gregorian calendar in which the  century years divisible by 400 are  not the leap years. 

$L_y$ is called \emph{lunar equation correction}, 
and it is equal to 1  only for years of the form $y=1800+ 300k+2500m$, where $k \in  \{0,1,2,\ldots,7\}$
and $m\in \{0,1,2,3\ldots\}$. 
This correction is necessary
to make the ecclesiastical moon as close as possible to the mean astronomical moon \cite{Roegel04themissing}.
The  sequence of
consecutive lunar equation correction years 
starts from 1800 and increases by 300 seven times, then a single increase by 400 follows, and afterwards the same sequence of increases is repeated \emph{ad infinitum}. We can, therefore, define
\begin{equation}\label{luncor}
 L_y = \begin{cases}
   1  & \text{if $\frac{(y-1800)\,\, \text{mod}\,\, 2500}{300}\in\{0,1,2,\ldots,7\}$,} \\
   0 &  \text{ otherwise.}
 \end{cases}
 \end{equation}
% \begin{equation}
% L_y = \begin{cases}
%   1  & \text{if $y$ mod $100=0$ and $(8C_{y}-3)$ mod $25 >16$,} \\
%   0 &  \text{ otherwise.}
% \end{cases}
% \end{equation}
% The rason why the above formula  works is because the remainder from
% division of $8n$ by $25$ for $n=1,2,\ldots, 25$ takes values in the set $$\{
% 8, 16, \underline{24}, 7, 15, \underline{23}, 6, 14, \underline{22}, 5, 13, \underline{21}, 4, 12, \underline{20}, 3, 11, \underline{19}, 2, 10, \underline{18}, 1, 9, \underline{17}, 0 \},$$ in which all integers from 0 to 24 occur once. Remainders greater than 16 are underlined, and we can see that there are 8 of them and they occur every thid position, exept that the last one (17) is separated by 4 positions form the next one (24) when the sequence is repeated.
% If instead of $8n \mod 25$ we take $8C_y-3 \mod 25$, a similar effect is obtained, namely $8C_y-3 \mod 25$ will be greater 
% than 16 for $C_n=19,22,25,28,31,34,37,40$,
% corresponding to years $1800, 2100, 2400, 2700, 3000, 3300, 3600, 3900$, as desired.
Recurrence equation (\ref{epactrecursion})  with
correction terms defined in eq. (\ref{metoncor}--\ref{luncor})
allows to compute the epact for any year after 1582.
Nevertheless, if the year is very large, the number of recursion steps needed is also very large, therefore we need
to find the solution of eq. (\ref{epactrecursion}), that is, 
the explicit formula for $E_y$ involving $y$.

\section{Solution of the recurrence equation for epacts}
The solution of the recurrence equation (\ref{epactrecursion}) can be constructed in two steps, starting from the simplified equation
\begin{align}
&E_{1582}=26, \nonumber \\
&E_y=E_{y-1}+11 + M_y \mod 30,
\end{align}
in which the solar and lunar corrections are omitted. 
Since we already know that the cycle of epacts returns to the same position after 19 years, to obtain $E_y$ we only need to add
to the initial value  26 the difference between golden numbers of $y$ and 
1582 multiplied by 11,
$$E_y=26+11 (G_y-G_{1582})\mod 30.$$
Since $G_{1582}=6$, this yields $11(G_y-G_{1582})+26=11G_y -40$,
and after adding 60 to ensure the expression is positive we obtain
\begin{equation}
E_y=11 G_y+20 \mod 30. 
\end{equation}
The above solution is actually also the correct solution 
of the full equation~(\ref{epactrecursion}) valid for years from 1582 to 1699 inclusively,
because in this period both solar and lunar corrections
remained zero.

In order to obtain the full solution for later years, we need to add the total sum of corrections,
\begin{equation}
 E_y=11 G_y+20 -\sum_{i=1583}^y S_i +\sum_{i=1583}^y L_i  \mod 30. 
\end{equation}
Let us start from the first sum. Since the solar correction can happen only in years divisible by 100, it will be convenient
to define century number $C_y$  of year $y$,
$$C_y=\lfloor y/100 \rfloor +1.$$
The total number of corrections from the year of the reform to $y$ can clearly only depend on $C_y$. Since $S_y$ is
equal to 1 three times in a period of 4 centuries,
if we were starting from year 0, the total number of corrections in the interval from 0 to $y$ inclusively would be $\lfloor 3C_y/4\rfloor$.  We are not starting from 0, however, so we need to subtract corrections in the period from 0 to 1852, thus
\begin{equation}
 \sum_{i=1583}^y S_i=\left\lfloor \frac{3 C_y}{4} 
 \right\rfloor -\left\lfloor \frac{3 C_{1582}}{4} 
 \right\rfloor.
\end{equation}
Since $C_{1582}=16$, the final result is
\begin{equation}
 \sum_{i=1583}^y S_i=\left\lfloor \frac{3 C_y}{4} 
 \right\rfloor -12.
\end{equation}
Now we will count the number of lunar corrections.
Consider at first the sequence $a_n=\lfloor 8n/25 \rfloor$ for $n=0,1,2, ..\ldots$. Its initial terms are
$$0, 0, 0, 0, 1, 1, 1, 2, 2, 2, 3, 3, 3, 4, 4, 4, 5, 5, 5, 6, 6, 6, 7, 7, 7, 8, 8, 8, 8,9\ldots$$
We can see that the sequence increases after 4 steps, then seven times after 3 steps, then again after 4 steps, etc.
The values of $a_n$ increase when $n=4+3m+25k$ where
$m\in \{0,1,\ldots,7\}$ and $k\in \{0,1,2, \ldots\}$.

This is exactly the kind of pattern we need to reproduce the sequence of  the years of the lunar correction.
Note that the first increase of $a_n$ happens at $n=4$, and the first lunar correction happens at $y=1800$, for which $C_n=19$,
meaning that $n=4$ corresponds to $C_n=19$. We can thus take
$n=C_y-15$, obtaining 
$$\left\lfloor \frac{ 8(C_y-15)}{25}\right\rfloor 
=\left\lfloor \frac{ 8C_y-120}{25}\right\rfloor 
=\left\lfloor \frac{ 8C_y+5-125}{25}\right\rfloor 
=\left\lfloor \frac{ 8C_y+5}{25}\right\rfloor -5.
$$
The final expression, therefore, is given by
\begin{equation}\label{totalL}
 \sum_{i=1583}^y L_i=\left\lfloor \frac{ 8C_y+5}{25}\right\rfloor -5.
\end{equation}
This yields the complete solution of eq. (\ref{epactrecursion}),
\begin{equation}\label{epactfinalraw}
E_y=11 G_y+20 -\left(\left\lfloor \frac{3 C_y}{4} 
 \right\rfloor -12\right) +
 \left\lfloor \frac{ 8C_y+5}{25}\right\rfloor -5 \mod 30,
\end{equation}
and, after simplification,
\begin{equation}\label{epactfinalsimplified}
E_y=11 G_y -\left\lfloor \frac{3 C_y}{4} 
 \right\rfloor  +
 \left\lfloor \frac{ 8C_y+5}{25}\right\rfloor +27 \mod 30.
\end{equation}

\section{Historical remarks}
\begin{figure}
\begin{center}
\includegraphics[width=14cm]{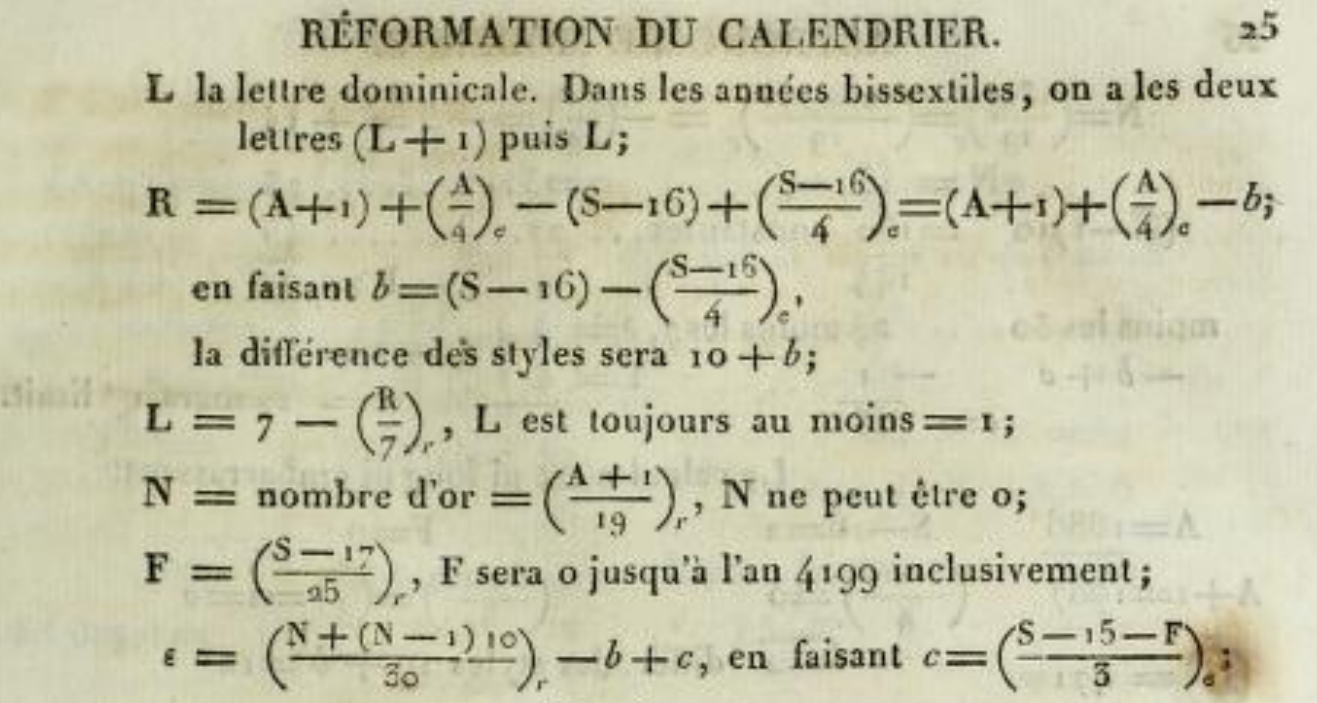}
\end{center}
\caption{Formulae for epacts's calculation from 
Delambre's \emph{Histoire de l'astronomie moderne} \cite{delambre1821}, p.25.
}\label{delambreformula}
\end{figure}

Before we continue, we will make some remarks about the
history of calculation of Gregorian epacts.
The ``canonical'' reference for the calculation of epacts (as well as the other aspects of the Gregorian calendar reform) is the work of Christopher Clavius \cite{Clavius1603}. 
Clavius described his computations using tables accompanied by extensive explanations of the method, but have not used any formulae. The first attempt to express the algorithm for epact calculation without any tables  is due to   Jean-Baptiste Delambre  \cite{delambre1814,delambre1821}. 
His formula for the Gregorian epact, taken
from \cite{delambre1821},  is shown in Figure
\ref{delambreformula}. It uses $A$ for year, $S$ for century
and $\epsilon$ for the epact. Symbol $(x/y)_e$ denotes the
integer part of $x/y$ and $(x/y)_r$ is the reminder from $x/y$ division.

\begin{figure}
\begin{center}
\includegraphics[width=12cm]{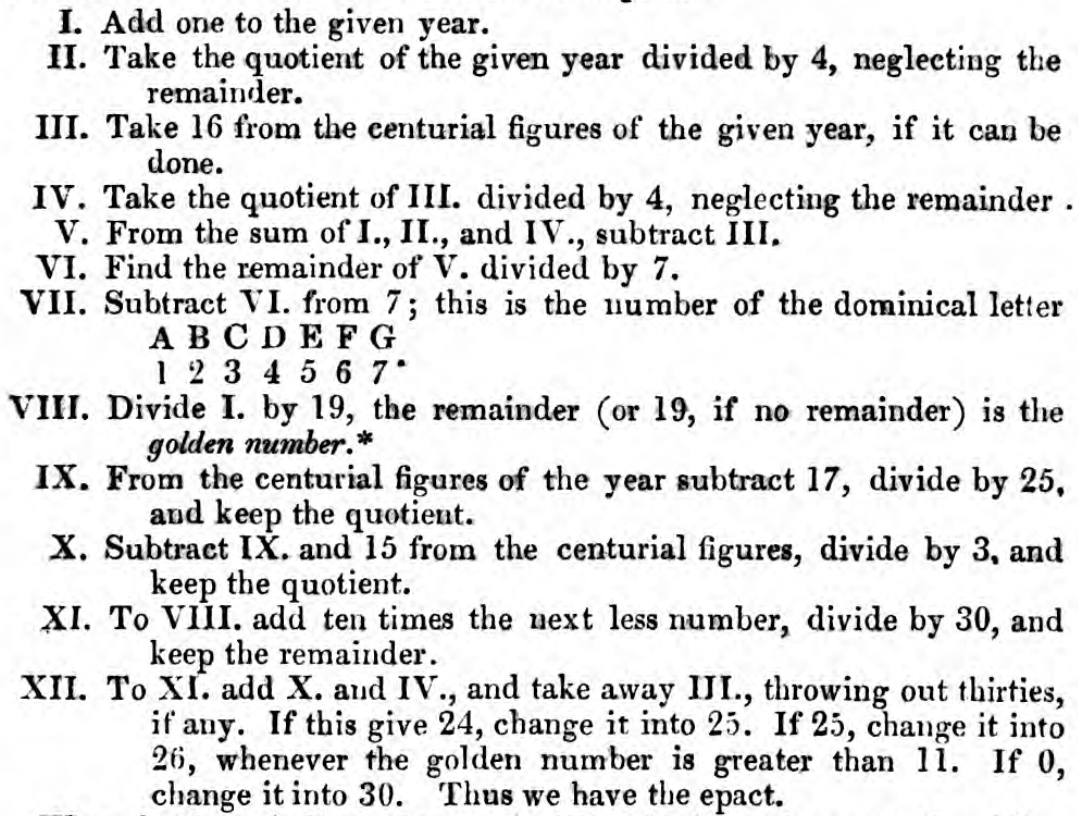}
\end{center}
\caption{De Morgan's Algorithm for epacts's calculation from 
\cite{Demorgan1845}. Only steps I and VIII-XII are relevant for the epact.
}\label{demorganalg}
\end{figure}

Delambre's formulae 
needed to calculate the epact of a given year has later been 
adopted by Augustus De Morgan in \cite{Demorgan1845}.
He took the formulae od Delambre 
and with some modifications he wrote them in an algorithmic form, as a sequence of 
artmetic operations, as shown in Figure \ref{demorganalg}.
This was a part of the algorithm for the computation of the date of Easter.

De Morgan later included his algorithm 
in his well known popularizing book \emph{Budget of paradoxes}
\cite{Demorgan1872}, and 
from there it made its way into numerous other
publications about the Gregorian calendar and the date of Easter, including the paper of  Donald Knuth~\cite{Knuth1962}. Knuth's formula for the epact is exactly as in De Morgan's work, and it is very similar to what we obtained here in
eq. (\ref{epactfinalraw}), except that the expression for the total lunar correction (which Knuth calls \emph{Clavian correction}) is somewhat different than our eq. (\ref{totalL}), namely it is
\begin{equation}\label{totLbig}
\sum_{i=1583}^y L_i=
\left\lfloor
\frac{C_n-\left\lfloor \frac{C_y-18}{25} \right\rfloor}{3}
\right\rfloor.
\end{equation}
One can easily show that this expression is equivalent to
eq. (\ref{totalL}) by simply checking all cases given by 
possible values of $C_y \mod 25$.
Interestingly, Knuth also gives the formula for the epact in~\cite{Knuth1997}, but there he already uses the form of the lunar correction exactly as in eq. (\ref{totalL}). There is no question that eq.  (\ref{totalL}), which uses properties of the sequence $\lfloor 8n/25 \rfloor$, is simpler and more elegant
than eq. (\ref{totLbig}), but who was first to use $\lfloor 8n/25 \rfloor$?  
It turns out that it is actually due to Gauss,
but you will not find it in his original paper on the algorithm for computing the date of Easter \cite{Gauss1800}. It appears
in the correction to that paper published in 1816 \cite{Gauss1816}, where Gauss gives the corrected form of 
one variables in his algorithm, $p=\lfloor (13+8k)/25 \rfloor$,
previously incorrectly given as $p=\lfloor k/3  \rfloor$.
Gauss is not calculating epacts explicitly, but he uses $p=\lfloor (13+8k)/25 \rfloor$ for the same 
purpose as we used  $\lfloor 8n/25 \rfloor$
, i.e., to account for solar equation corrections.
Afterwards the use of the sequence $\lfloor 8n/25 \rfloor$ to calculate the lunar correction became widely know and it has been adopted to other algorithms, including some recent 
versions, such as \cite{Lichtenberg1997}.

\section{Age of the  ecclesiastical moon}
Having the epact, the computation of the age of the eccelsiastical moon is straightforward.
The epact of the year $E_y$ plus one gives the age of the moon on  January 1. The following principles, which can be deduced from the \emph{Calendarium} and from the lunar tables of the \emph{Martyrologium}, are used to calculate the age of the moon at other  days of the year.
\begin{itemize}
 \item the age of the moon in the first days of January increases till it reaches 30, and then it  follows one of the two sequences, either
 $$a=\{1,2,\ldots,29, 1,2, \ldots, 30, 1,2,\ldots,29,\ldots,\}$$
 or
 $$b=\{1,2,\ldots,30, 1,2, \ldots, 29, 1,2,\ldots,30,\ldots,\}$$
 that is, alternating sequence of cycles 29 and 30 days long.
 \item If $E_y<25$ or ($E_y=25$ and $G_y\geq 12$) then the first of the above sequences is followed.
 \item in all other cases, the second  sequence is followed.
\end{itemize}
In a year with epact $E_y$ the consecutive ages of the moon will then be
$$
 E_y+1,E_y+2,\ldots,30, 1,2,\ldots,29, 1,2, \ldots, 30, 1,2,\ldots,29,\ldots  
$$
or
$$
 E_y+1,E_y+2,\ldots,30, 1,2,\ldots,30, 1,2, \ldots, 29, 1,2,\ldots,29,\ldots  
$$
The two sequences were introduced by Clavius in the design of the calendar in order to avoid certain regularities of the ecclesiastical moon, namely to avoid a possibility of having new moon on the same date in two different years within the same 
Metonic cycle as well as avoiding Easter Sunday to occur on the same day of the month twice in 19 years (see \cite{Demorgan1845} for more details and explanation). What is important is that both sequences will result in the same 
average synodic month duration of $29.5$ days, as desired.

Formula for the age of the moon can be now constructed.
Let $g(x)$ be a function such that its values for $x=0,1,2,3\ldots$ are exactly as in the sequence  $a$, so that
$g(0)=1$, $g(1)=2$, \ldots $g(28)=29$, $g(30)=1$, $g(31)=2$, and then it cycles to 29, then to 30, again to 29 etc.
It is easy to see that for some fixed $k$,
when $x$ take consecutive values $0,1,2,$ the values of
$g(x+k)$ will be follow the pattern
$$
 k+1,k+2,\ldots,30, 1,2,\ldots,29, 1,2, \ldots, 30, 1,2,\ldots,29,\ldots  
$$
This is exactly what we need for the age of the moon following case $a$ if we take $k=E_y$.

For sequence $b$ we first note that we can obtain $b$ by shifting $a$ by $29$. Furthermore, for fixed $k$,
when $x$ take consecutive values $0,1,2,$ values of
$g(x+k+29)$ will be follow the pattern
$$
 k,k+1,k+2,\ldots,29, 1,2,\ldots,30, 1,2, \ldots, 29, 1,2,\ldots,29,\ldots  
$$
This is almost the desired sequence of moon ages following pattern $b$, except that we need to add 1 to the initial part
$k, k+1,\ldots, 29$ in order to make it $k+1, k+2,\ldots, 30$,
and then again take $k=E_y$.

We can now summarize our findings as follows.
Let $n$ be the consecutive number of the day of the year starting from $n=0$ on January 1, $n=1$ on January 2 etc.
February 29 in leap years is omitted, and we assume that
February 29 has the same day number as February 28.
The age of the moon $A_{y,n}$ at year $y$  and day $n$ is then
\begin{equation}\label{moonage}
A_{y,n} = \begin{cases}
  g(E_y+n)  & \text{if $E_y<25$ or ($E_y=25$ and $G_y\geq 12$),} \\
 g(E_y+n+29)  + \mathbf{1}_{n+E_y<30}
 &  \text{ otherwise,}
\end{cases}
\end{equation}
where
\begin{equation}
\mathbf{1}_{P}=
\begin{cases}
 1  & \text{if $P$ is true} \\
 0  &  \text{ otherwise.}
\end{cases}
\end{equation}
All what is left is to find a formula for $g(x)$.
Here is one possibility,
\begin{equation}
g(x)=\left(x+ \left\lfloor \frac{x}{59}  \right\rfloor \right) \text{\,mod\,} 30  \,\,\, + 1.
\end{equation}
The reader can easily verify that $g$ is periodic with period 59, and that within a single period $g(x)$ increases
from 1 to 29 and then from 1 to 30, just as required. 

To avoid any confusion, let us also provide the formula for
the day number in terms of month  $m$ and day $d$.
The day number, previously denoted by $n$ but to be further denoted by $N_{m,d}$, can be
calculated as follows:
\begin{equation}
 N_{m,d}=d-1 + 30(m-1)+\left\lfloor \frac{7m-2}{12}\right\rfloor
-2 \left\lfloor \frac{m+9}{12}\right\rfloor.
\end{equation}
In order to understand where this is coming from, let us note that
in the above, $d-1 + 30(m-1)$ would be sufficient if all months had the length 30 days. The first correction,
$\left\lfloor \frac{7m-2}{12}\right\rfloor$, takes care of the months which have the length 31 days, and the last term,
$-2 \left\lfloor \frac{m+9}{12}\right\rfloor$ corrects
for February having 28 days. It should be stressed that the above formula neglects the difference between common and leap years, and in the leap years we simply assume that $N_{2,29}=N_{2,28}$.

The final formula for the age of the ecclesiastical moon $A_{y,m,d}$  on the day $d$ of  month $m$ and year $y$ is given by
\begin{equation} \label{ageofmoon}
A_{y,m,d} = \begin{cases}
  g(E_y+N_{m,d})  & \text{if $E_y<25$ or ($E_y=25$ and $G_y\geq 12$),}  \\
 g(E_y+N_{m,d}+29)  + \mathbf{1}_{N_{m,d}+E_y<30}
 &  \text{ otherwise.}
\end{cases}
\end{equation}

\section{Corrections at the beginning and end of the year}
Let us now consider what happens at the transition from December 31 of year $y-1$ to January 1 of year $y$.
It is easy to check that on December 31 of year $y-1$ the age of the moon
is always
$$
A_{y-1,12,31}=(E_{y-1}+10) \mod 30 +1.
$$
This can be verified by direct substitution of $N_{12,31}=364$ into
eq.~(\ref{ageofmoon}). By similar substitution of $N_{1,1}=0$ 
into the same equation we obtain the age of the moon on January 1 of the year $y$
$$
A_{y,1,1}=E_{y}+1.
$$
This is, of course, not surprising because this is how the epact has been defined (age of the moon on January 1 minus 1).
According to eq. (\ref{epactrecursion})
$$
E_{y}=E_{y-1}+11 + M_{y}- S_{y}+L_{y} \mod 30,
$$
therefore, if $M_{y}= S_{y}=L_{y}=0$, we obtain
$$
A_{y,1,1}=(E_{y-1}+11) \mod 30 +1,
$$
meaning that the moon is one day older on Jan 1 than it was on Dec. 31 of the previous year. This is how it should be, but this will 
unfortunately not work when  $M_{y}- S_{y}+L_{y} \neq 0$. In such a case there might be some jump in the age of the moon, by a day or even two. The liturgical books deal with this problem partially, by providing a correction which works for $M_{y}=1$, $- S_{y}+L_{y}=0$ case, that is, for the transition from the year with golden number 
19 to they year with $G_y=1$, when solar and lunar corrections are both zero or cancel each other.

Curiously, there are two fixes appearing in two different places, and these are mutually exclusive, meaning that they should not be applied together.

The first appears in the \emph{Calendarium} section of the Breviary (and also Missal). In the column where we find epacts for December 31 we find that there two epacts, one being the usual xx, and the other 19, printed in Arabic numerals. The explanation below the table says that the epact 19 is to be used only when the golden number of the year is also 19.

Why it fixes the problem is best explained on example. In 1690  the
epact was xix and the golden number was 19. Epact xix appears in \emph{Calendarium}
on December 2, which is therefore the first day of the moon. December 31, if we count down, is then 30-th day of the moon.

Now check the next year, 1691, which has golden number 1 and epact i.
This means that on January 1 the age of the moon is 2. The new moon (of age 1) is, therefore, missing. In order to fix this \emph{Calendarium}
requests to place an additional epact 19 on Dec. 31, meaning that this
will be the day of the new moon (of age 1). December 30 will retain
age of the moon 29, so this simply shortens the December cycle to 29 days.

In the Martyrology the problem is fixed differently, by shifting the moon age in January of the next year by one day. The introductory section of the Martyrology advises that for years of the golden number 1, excluding years with epacta *, the age of the of the moon from January 1 to the end of the cycle should be ``pronounced'' as decreased by 1 day.
In our example year of 1691 this means that ages of the moon 
in January will be pronounced $1,2,3,4,\ldots$ instead of $2,3,4,\ldots$, obtained from the lunar tables of the Martyrologium. This way
the missing new moon is restored on January 1. A formal
way of writing this corrected moon age, denoted by $A_{y,m,d}^{pron}$, would be
\begin{equation*}
A_{y,m,d}^{pron} = \begin{cases}
 A_{y,m,d}-1 &\text{if $G_y=1$, $E_y>0$, $m=1$ and  $d+E_y \leq 30$,} \\
 A_{y,m,d} &  \text{ otherwise.}
\end{cases}
\end{equation*}
In the above the condition $d+E_y \leq 30$ appears because we want the
correction to be applied from $d=1$ to $d=30-E_y$ (end of the cycle).
Condition $E_y>0$ appears because when $E_y=0$, the moon age on
January 1 is 1, so no correction is needed.

To demonstrate how this works consider the next time we will have $G_{y}=1$, which is 2033.  The table of moon ages for December of 2032 and January 2033, constructed using uncorrected $A_{y,m,d}$,
will then look like this (with new moon shown in red and 14th day in blue):
\begin{center}
\small
\setlength{\tabcolsep}{2.3pt}
\renewcommand{\arraystretch}{1.2}
\begin{tabular}{c|ccccccccccccccccccccccccccccccc}
   & 1& 2& 3& 4& 5& 6& 7& 8& 9&10&11&12&13&14&15&16&17&18&19&20&21&22&23&24&25&26&27&28&29&30&31\\
\hline
Dec  &27&28&29&{\color{red} 1}& 2& 3& 4& 5& 6& 7& 8& 9&10&11&12&13&{\color{blue}14}&15&16&17&18&19&20&21&22&23&24&25&26&27&28\\
Jan  &30&{\color{red} 1}& 2& 3& 4& 5& 6& 7& 8& 9&10&11&12&13&{\color{blue}14}&15&16&17&18&19&20&21&22&23&24&25&26&27&28&29&30
\end{tabular}
\end{center}
We can see that the moon age on Dec. 31 is 28, and it is followed by 30 on January 1. There is a jump from 28 to 30, and day 29 is lost. When we apply the Martyrology correction  and use $A^{corr}_{y,m,d}$ to generate the same table, the result for  2032/2033 will look like this:
\begin{center}
\small
\setlength{\tabcolsep}{2.3pt}
\renewcommand{\arraystretch}{1.2}
\begin{tabular}{c|ccccccccccccccccccccccccccccccc}
   & 1& 2& 3& 4& 5& 6& 7& 8& 9&10&11&12&13&14&15&16&17&18&19&20&21&22&23&24&25&26&27&28&29&30&31\\
\hline
Dec  &27&28&29&{\color{red} 1}& 2& 3& 4& 5& 6& 7& 8& 9&10&11&12&13&{\color{blue}14}&15&16&17&18&19&20&21&22&23&24&25&26&27&28\\
Jan  &29&{\color{red} 1}& 2& 3& 4& 5& 6& 7& 8& 9&10&11&12&13&{\color{blue}14}&15&16&17&18&19&20&21&22&23&24&25&26&27&28&29&30
\end{tabular}
\end{center}

Note that when we use the Martyrology correction, we shall
not use the solution from the Missal at the same time, 
because if we did, we would have two new moons on consecutive days, on Dec. 31 and Jan 1.
To see this, consider the next case of the year when $G_y=E_y=19$, $y=8511$.
If we generate the table of moon ages using uncorrected $A_{y,m,d}$ for 
 8511/8512, we obtain
\begin{center}
\small
\setlength{\tabcolsep}{2.3pt}
\renewcommand{\arraystretch}{1.2}
\begin{tabular}{c|ccccccccccccccccccccccccccccccc}
   & 1& 2& 3& 4& 5& 6& 7& 8& 9&10&11&12&13&14&15&16&17&18&19&20&21&22&23&24&25&26&27&28&29&30&31\\
\hline
Dec  &29&{\color{red} 1}& 2& 3& 4& 5& 6& 7& 8& 9&10&11&12&13&{\color{blue}14}&15&16&17&18&19&20&21&22&23&24&25&26&27&28&29&30\\
Jan  & 2& 3& 4& 5& 6& 7& 8& 9&10&11&12&13&{\color{blue}14}&15&16&17&18&19&20&21&22&23&24&25&26&27&28&29&30&{\color{red} 1}& 2
\end{tabular}
\end{center}
The same table for 8511/8512, generated using the pronounced moon age $A^{pron}_{y,m,d}$ will become
\begin{center}
\small
\setlength{\tabcolsep}{2.3pt}
\renewcommand{\arraystretch}{1.2}
\begin{tabular}{c|ccccccccccccccccccccccccccccccc}
   & 1& 2& 3& 4& 5& 6& 7& 8& 9&10&11&12&13&14&15&16&17&18&19&20&21&22&23&24&25&26&27&28&29&30&31\\
\hline
Dec  &29&{\color{red} 1}& 2& 3& 4& 5& 6& 7& 8& 9&10&11&12&13&{\color{blue}14}&15&16&17&18&19&20&21&22&23&24&25&26&27&28&29&30\\
Jan  &{\color{red} 1}& 2& 3& 4& 5& 6& 7& 8& 9&10&11&12&13&{\color{blue}14}&15&16&17&18&19&20&21&22&23&24&25&26&27&28&29&{\color{red} 1}& 2
\end{tabular}
\end{center}
If we used the ``black epact 19'' fix  proposed in the \emph{Calendarium} and applied it  to the above, we would have to change the moon age on Dec 31 from 30 to 1, obtaining two consecutive new moons, as shown below.
\begin{center}
\small
\setlength{\tabcolsep}{2.3pt}
\renewcommand{\arraystretch}{1.2}
\begin{tabular}{c|ccccccccccccccccccccccccccccccc}
   & 1& 2& 3& 4& 5& 6& 7& 8& 9&10&11&12&13&14&15&16&17&18&19&20&21&22&23&24&25&26&27&28&29&30& 31\\
\hline
Dec  &29&{\color{red} 1}& 2& 3& 4& 5& 6& 7& 8& 9&10&11&12&13&{\color{blue}14}&15&16&17&18&19&20&21&22&23&24&25&26&27&28&29&{\color{red}1}\\
Jan  &{\color{red} 1}& 2& 3& 4& 5& 6& 7& 8& 9&10&11&12&13&{\color{blue}14}&15&16&17&18&19&20&21&22&23&24&25&26&27&28&29&{\color{red} 1}& 2
\end{tabular}
\end{center}

As mentioned, the ``pronounced'' moon age, while it solves the problem partially, works only when there is a jump by one day at the transition from
a year with $G_{y-1}=19$ to $G_y=1$. It does not solve the problem of
other anomalies, for example, jump by two days when both $M_{y}=S_{y}=1$ and $L_y=0$, or ``stalling'' of the moon when  $M_{y}=S_{y}=0$ and $L_y=1$. A discussion of these anomalies can be found in \cite{Roegel04themissing}. The full fix for all these cases is proposed in the next section.

\section{Proposed correction which takes care of all cases}
It is straightforward to construct a corrected moon age
following the idea of Martyrologium's correction. 
Looking at eq.~(\ref{epactrecursion}) we can see that the 
``jump'' between Dec. 31 and Jan 1  is equal to  $M_y-Sy+Ly$
and can take values  -1, 0, 1 or 2, with 0 being the regular
case (no jump). Therefore, in order to eliminate the jump, 
we simply need to subtract the value of the jump from
moon ages for all days from January 1 to the end of the cycle, 
that is, for $m=1$ and $d=1, \ldots 30-E_y$.

The first thing we need is to calculate the value of the jump $J_y =M_y-Sy+Ly$, and the most 
convenient way to do this is to use our formula for the epact $E_y$.
From eq.~(\ref{epactrecursion}), the jump is just the difference between $E_y$ and $E_y-1$, minus 11. We just need to take care of the fact that epacts are taken modulo 30, thus the correct formula for the jump will be
$$
J_y=(E_y-E_{y-1}) \text{\,\,\,mod 30}\,\,\, -11.
$$
Now all what is needed is decreasing $A_{y,m,d}$ by $J_y$ on 
on days from Jan. 1 to the end of the cycle. However, we
cannot just write $A_{y,m,d}- J_y$, as it may happen that
$A_{y,m,d}\leq J_y$, and then we would get zero or negative day of the moon. For this reason, we need to bring the
result to the range $1\ldots 30$ by adding 30,
which yields the formula for the corrected moon age as
\begin{equation*}
A_{y,m,d}^{corr} = \begin{cases}
%(ageofmoon(y,m,d)-jump -1)%30 +1
%jump!=0 and  E>0 and m==1 and d<=30-E :
A_{y,m,d}-J_y +30\cdot \mathbf{1}_{A_{y,m,d}-J_y \leq 0} &\text{if  $m=1$ and  $d+E_y \leq 30$,} \\
 A_{y,m,d} &  \text{ otherwise.}
\end{cases}
\end{equation*}

This correction fixes all anomalies listed in \cite{Roegel04themissing}. Let us discuss some examples.

The first one will be the case of the year 16399/16400, mentioned in the title of \cite{Roegel04themissing}.
For $y=16400$ we have $E_{y-1}=19$, $E_{y}=1$, $J_{y}=1$. We have a jump by one day, and since $G_y=4$, it
would not be fixed by $A^{pron}_{y,m,d}$. Uncorrected table of moon ages would look like this:
\begin{center}
\small
\setlength{\tabcolsep}{2.3pt}
\renewcommand{\arraystretch}{1.2}
\begin{tabular}{c|ccccccccccccccccccccccccccccccc}
   & 1& 2& 3& 4& 5& 6& 7& 8& 9&10&11&12&13&14&15&16&17&18&19&20&21&22&23&24&25&26&27&28&29&30&31\\
\hline
Dec  &29&{\color{red} 1}& 2& 3& 4& 5& 6& 7& 8& 9&10&11&12&13&{\color{blue}14}&15&16&17&18&19&20&21&22&23&24&25&26&27&28&29&30\\
Jan  & 2& 3& 4& 5& 6& 7& 8& 9&10&11&12&13&{\color{blue}14}&15&16&17&18&19&20&21&22&23&24&25&26&27&28&29&30&{\color{red} 1}& 2
\end{tabular}
\end{center}
The new moon is missing on Jan 1. If we apply the proposed correction and generate
the same table for 16399/16400 using $A^{corr}_{y,n,d}$, the problem will be fixed:
\begin{center}
\small
\setlength{\tabcolsep}{2.3pt}
\renewcommand{\arraystretch}{1.2}
\begin{tabular}{c|ccccccccccccccccccccccccccccccc}
   & 1& 2& 3& 4& 5& 6& 7& 8& 9&10&11&12&13&14&15&16&17&18&19&20&21&22&23&24&25&26&27&28&29&30&31\\
\hline
Dec  &29&{\color{red} 1}& 2& 3& 4& 5& 6& 7& 8& 9&10&11&12&13&{\color{blue}14}&15&16&17&18&19&20&21&22&23&24&25&26&27&28&29&30\\
Jan  &{\color{red} 1}& 2& 3& 4& 5& 6& 7& 8& 9&10&11&12&13&{\color{blue}14}&15&16&17&18&19&20&21&22&23&24&25&26&27&28&29&{\color{red} 1}& 2
\end{tabular}
\end{center}

Another example is the year 106399/106400, for which $E_{y-1}=18$, $E_{y}=1$, $J_{y}=2$. Here is the uncorrected table:
\begin{center}
\small
\setlength{\tabcolsep}{2.3pt}
\renewcommand{\arraystretch}{1.2}
\begin{tabular}{c|ccccccccccccccccccccccccccccccc}
   & 1& 2& 3& 4& 5& 6& 7& 8& 9&10&11&12&13&14&15&16&17&18&19&20&21&22&23&24&25&26&27&28&29&30&31\\
\hline
Dec  &28&29&{\color{red} 1}& 2& 3& 4& 5& 6& 7& 8& 9&10&11&12&13&{\color{blue}14}&15&16&17&18&19&20&21&22&23&24&25&26&27&28&29\\
Jan  & 2& 3& 4& 5& 6& 7& 8& 9&10&11&12&13&{\color{blue}14}&15&16&17&18&19&20&21&22&23&24&25&26&27&28&29&30&{\color{red} 1}& 2
\end{tabular}
\end{center}
As we can see, not only there is a jump by two days, but the new moon (of age 1) is missing at the beginning of January 106400. The same table generated using $A_{y,m,d}^{corr}$ fixes this fully:
\begin{center}
\small
\setlength{\tabcolsep}{2.3pt}
\renewcommand{\arraystretch}{1.2}
\begin{tabular}{c|ccccccccccccccccccccccccccccccc}
   & 1& 2& 3& 4& 5& 6& 7& 8& 9&10&11&12&13&14&15&16&17&18&19&20&21&22&23&24&25&26&27&28&29&30&31\\
\hline
Dec  &28&29&{\color{red} 1}& 2& 3& 4& 5& 6& 7& 8& 9&10&11&12&13&{\color{blue}14}&15&16&17&18&19&20&21&22&23&24&25&26&27&28&29\\
Jan  &30&{\color{red} 1}& 2& 3& 4& 5& 6& 7& 8& 9&10&11&12&13&{\color{blue}14}&15&16&17&18&19&20&21&22&23&24&25&26&27&28&{\color{red} 1}& 2
\end{tabular}
\end{center}

The third and the last example we will show is for year 4199/4200, with $E_{y-1}=20$, $E_{y}=0$, $J_{y}=-1$, when we have two consecutive new moon days on Dec. 31 and Jan 1 in the uncorrected table:
\begin{center}
\small
\setlength{\tabcolsep}{2.3pt}
\renewcommand{\arraystretch}{1.2}
\begin{tabular}{c|ccccccccccccccccccccccccccccccc}
   & 1& 2& 3& 4& 5& 6& 7& 8& 9&10&11&12&13&14&15&16&17&18&19&20&21&22&23&24&25&26&27&28&29&30&31\\
\hline
Dec  &{\color{red} 1}& 2& 3& 4& 5& 6& 7& 8& 9&10&11&12&13&{\color{blue}14}&15&16&17&18&19&20&21&22&23&24&25&26&27&28&29&30&{\color{red} 1}\\
Jan  &{\color{red} 1}& 2& 3& 4& 5& 6& 7& 8& 9&10&11&12&13&{\color{blue}14}&15&16&17&18&19&20&21&22&23&24&25&26&27&28&29&30&{\color{red} 1}
\end{tabular}
\end{center}
Again, $A_{y,m,d}^{corr}$ fixes this problem, and the corrected table for 4199/4200 is as follows:
\begin{center}
\small
\setlength{\tabcolsep}{2.3pt}
\renewcommand{\arraystretch}{1.2}
\begin{tabular}{c|ccccccccccccccccccccccccccccccc}
   & 1& 2& 3& 4& 5& 6& 7& 8& 9&10&11&12&13&14&15&16&17&18&19&20&21&22&23&24&25&26&27&28&29&30&31\\
\hline
Dec  &{\color{red} 1}& 2& 3& 4& 5& 6& 7& 8& 9&10&11&12&13&{\color{blue}14}&15&16&17&18&19&20&21&22&23&24&25&26&27&28&29&30&{\color{red} 1}\\
Jan  & 2& 3& 4& 5& 6& 7& 8& 9&10&11&12&13&{\color{blue}14}&15&16&17&18&19&20&21&22&23&24&25&26&27&28&29&30&31&{\color{red} 1}
\end{tabular}
\end{center}
Note that $A^{corr}_{y,m,d}$, compared to $A_{y,m,d}$, changes
only the first lunation of January on years of $J_y \neq 0$, yet
it leaves all other lunations intact. For this reason, it will not affect the date of Easter or dates or other movable feasts.
\section{Summary}
Our results can be summarized as follows.
For a given day $d$ of month $m$ and year $y$, define, respectively,
the century number $C_y$, the golden number $G_y$, the epact $E_y$ and the day number $N_{m,d}$ as follows,
\begin{align*}
 C_y&=\lfloor y/100 \rfloor +1,\\
 G_y&=y \mod 19 +1,\\
 E_y&=11 G_y -\left\lfloor \frac{3 C_y}{4} 
 \right\rfloor  +
 \left\lfloor \frac{ 8C_y+5}{25}\right\rfloor +27 \mod 30,\\
  N_{m,d}&=d-1 + 30(m-1)+\left\lfloor \frac{7m-2}{12}\right\rfloor
-2 \left\lfloor \frac{m+9}{12}\right\rfloor.
\end{align*}
Then the age of the ecclesiastical moon $A_{y,m,d}$  on that day is given by
\begin{equation*}
A_{y,m,d} = \begin{cases}
  g(E_y+N_{m,d})  & \text{if $E_y<25$ or ($E_y=25$ and $G_y\geq 12$),} \\
 g(E_y+N_{m,d}+29)  + \mathbf{1}_{N_{m,d}+E_y<30}
 &  \text{ otherwise,}
\end{cases}
\end{equation*}
% \begin{equation*}
% A_{y,m,d} = \begin{cases}
%   g(E_y+N_{m,d})  & \text{if $E_y<25$ or ($E_y=25$ and $G_y\geq 12$),} \\
%  g(E_y+N_{m,d}+29) 
%  + 1-\left\lfloor \frac{337+N_{m,d}+E_y}{366} \right\rfloor
%  &  \text{ otherwise,}
% \end{cases}
% \end{equation*}
where 
\begin{equation*}
g(x)=\left(x+ \left\lfloor \frac{x}{59}  \right\rfloor \right) \text{\,mod\,} 30  \,\,\, + 1.
\end{equation*}
The ``pronounced'' age of the moon, following \emph{Martyrologium Romanum}, is given by
\begin{equation*}
A_{y,m,d}^{pron} = \begin{cases}
 A_{y,m,d}-1 &\text{if $G_y=1$, $E_y>0$, $m=1$ and  $d+E_y \leq 30$,} \\
 A_{y,m,d} &  \text{ otherwise.}
\end{cases}
\end{equation*}
The proposed fully corrected moon age is given by
\begin{align*}
J_y&=(E_y-E_{y-1}) \text{\,\,\,mod 30}\,\,\, -11,\\
A_{y,m,d}^{corr} &= \begin{cases}
%(ageofmoon(y,m,d)-jump -1)%30 +1
%jump!=0 and  E>0 and m==1 and d<=30-E :
A_{y,m,d}-J_y +30\cdot \mathbf{1}_{A_{y,m,d}-J_y \leq 0} &\text{if  $m=1$ and  $d+E_y \leq 30$,} \\
 A_{y,m,d} &  \text{ otherwise.}
\end{cases}
\end{align*}

The listing below shows the implementation in Python of functions 
for the epact, the moon age, the pronounced moon age and the proposed fully corrected moon age.
We hope that the proposed correction
makes its way into the future edition of \emph{Martyrologium  Romanum} well before the year
16399.

\vskip 0.5cm
\begin{python}
 def epacta(y):
    C= y//100 + 1 	 
    G= y % 19 + 1 	
    E= (11*G -(3*C)//4 + (8*C+5)//25 + 27)%30 
    return E   
    
 def ageofmoon(y,m,d):
    C= y//100 + 1 	 
    G= y % 19 + 1 	 
    E= (11*G -(3*C)//4 + (8*C+5)//25 + 27)%30 
    N=d-1+30*(m-1)+ (7*m-2)//12 -2*((m+9)//12) 
    if E<25 or (E==25 and G>=12):
        return 1 + (N+E + (N+E)//59)%30
    return 1 + (N+E+29+ (N+E+29)//59)%30  + (N<30-E)
    
 def ageofmoon_pronounced(y,m,d):
     G=y % 19 + 1
     C= y//100 + 1 	 
     E= (11*G -(3*C)//4 + (8*C+5)//25 + 27)%30 
     if G==1 and  E>0 and m==1 and d<=30-E : 
           return ageofmoon(y,m,d)-1
     return ageofmoon(y,m,d)  
    
  def ageofmoon_fullycorrected(y,m,d):
     E=epacta(y)
     Eprev=epacta(y-1)
     J=(E-Eprev)%30 -11
     A=ageofmoon(y,m,d)
     if J!=0 and m==1 and d<=30-E : 
           return A-J +30*((A-J)<=0)
     return A  
\end{python}
\

\end{document}